# Exergy analysis of marine waste heat recovery $CO_2$ closed-cycle gas turbine system

Vedran Mrzljak[1], Igor Poljak[2], Jasna Prpić-Oršić[1], Maro Jelić[3]

[1]Faculty of Engineering, University of Rijeka, Vukovarska 58, 51000 Rijeka, Croatia
[2]Department of Maritime Sciences, University of Zadar, Mihovila Pavlinovića 1, 23000 Zadar, Croatia
[3]Maritime Department, University of Dubrovnik, Ćira Carića 4, 20000 Dubrovnik, Croatia
email: vedran.mrzljak@riteh.hr, ipoljak1@unizd.hr, jasna.prpic-orsic@riteh.hr, maro.jelic@unidu.hr

**Abstract:** This paper presents an exergy analysis of marine waste heat recovery $CO_2$ closed-cycle gas turbine system. Based on the operating parameters obtained in system exploitation, it is performed analysis of each system component individually, as well as analysis of the whole observed system. While observing all heat exchangers it is found that combustion gases-$CO_2$ heat exchangers have the lowest exergy destructions and the highest exergy efficiencies (higher than 92%). The lowest exergy efficiency of all heat exchangers is detected in Cooler (51.84%). Observed system is composed of two gas turbines and two compressors. The analysis allows detection of dominant mechanical power producer and the dominant mechanical power consumer. It is also found that the turbines from the observed system have much higher exergy efficiencies in comparison to compressors (exergy efficiency of both turbines is higher than 94%, while exergy efficiency of both compressors did not exceed 87%). The whole observed waste heat recovery system has exergy destruction equal to 6270.73 kW, while the exergy efficiency of the whole system is equal to 64.12% at the selected ambient state. Useful mechanical power produced by the whole system and used for electrical generator drive equals 11204.80 kW. The obtained high exergy efficiency of the whole observed system proves its application on-board ships.

## 1 Introduction

Marine propulsion systems are nowadays mainly based on the internal combustion diesel engines [1, 2]. Due to its dominancy, scientists and researchers are developing various simulation models which can accurately and precisely track diesel engines operating parameters [3, 4] with an aim to improve its operation. One of the dominant topics in the field of marine internal combustion diesel engines is reducing of its emissions, therefore various systems and methods are developed to satisfy harmful emissions legislative [5-7].

Other marine propulsion systems are significantly less represented in the entire worldwide fleet. As for example, steam propulsion systems are still dominantly represented in the propulsion of LNG (Liquefied Natural Gas) carriers due to the specificity of its operations and transported cargo [8, 9]. However, internal combustion dual-fuel engines each day takes stronger and stronger impact also in LNG carriers propulsion systems, therefore it can be expected that in the near future it will overcome traditional steam propulsion [10, 11].

Gas turbines as a stand-alone devices are rarely used for marine propulsion due to its high price, complex maintenance and other elements which can cause problems in operation or for a ship crew [12]. Gas turbines in marine propulsion can be in the most of the cases found as part of complex propulsion systems which includes several different propulsion elements [13, 14].

The second option for application of the gas turbine in the marine propulsion is using a various upgrades and modifications which increase the efficiency of the whole system [15]. Such upgrades can be performed on the gas turbine [16] or can be installed as an additional component [17, 18]. Many techniques are developed to ensure application of heat from gas turbine combustion gases for various heating purposes or for additional electrical power production [19-21]. Waste heat recovery systems allow that the whole propulsion plant with a gas turbine as a main propulsion element become a valid competitor to all other propulsion systems [22, 23].

In this paper is performed exergy (second law analysis) of $CO_2$ closed-cycle gas turbine system which for its operation uses heat from gas turbine combustion gases. This system is a good representation of the widely used techniques for waste heat recovery in additional electrical power production. Based on the obtained operating parameters from exploitation, it is performed calculation of exergy power inputs, outputs, exergy destructions (exergy losses) and exergy efficiencies of each system component and of the whole observed system. By taking into consideration many challenges which will surely occur in practical implementation of the observed system, obtained exergy efficiency of the analyzed waste heat recovery system makes the whole marine propulsion plant (gas turbine + presented $CO_2$ closed-cycle gas turbine system) a good competitor in comparison to other marine propulsion systems.

## 2 Description and operating characteristics of the analyzed $CO_2$ closed-cycle gas turbine system

Scheme of the analyzed waste heat recovery $CO_2$ closed-cycle gas turbine system is presented in Fig. 1. In Fig. 1 are also presented all the necessary operating points required for the exergy analysis of each system component. Red lines represent the flow of combustion gases, light blue dashed lines represent the flow of $CO_2$, while green dotted lines represent the flow of cooling water in the Cooler.

Analyzed system is composed of several heat exchangers [24, 25], two $CO_2$ turbines and two $CO_2$ compressors. The dominant heat exchangers are H1 and H2, in which combustion gases from the main marine gas turbine (and if required form the additional combustion chamber) is used for $CO_2$ heating before its expansion in T1 and T2. HTR, LTR and IHX are $CO_2$-$CO_2$ heat exchangers in which $CO_2$ of a higher temperature transfer heat to $CO_2$ of a lower temperature. Two compressors (Compr1 and Compr2) are used for the $CO_2$ pressure increase. $CO_2$ cooling is performed in only one component (Cooler) and for the $CO_2$ cooling purposes is used water ($CO_2$ transfer heat to cooling water), Fig. 1.

Both turbines and compressors are mounted on the same shaft. A part of cumulative mechanical power produced in both gas turbines is firstly used for the drive of both compressors, while the remaining part of cumulative produced mechanical power (useful power) is used for the electrical generator drive.

In a standard operation, the system did not require any additional fuel consumption – additional combustion chamber is not in operation, both heaters H1 and H2 uses combustion gases produced in a gas turbine only. At the highest loads, when the heat in combustion gases from the gas turbine is not sufficient for required heating in H1 and H2, the additional combustion chamber is applied.

It should be highlighted that in Fig. 1 are not shown $CO_2$ tanks which are used for the change of $CO_2$ mass flow rate through the observed system and compressors for $CO_2$ delivery from the observed system to that tanks or in a reverse direction. By changing of $CO_2$ mass flow rate through the observed system can be performed regulation of produced mechanical power (the same technique is valid for any other closed-cycle gas turbine [26]).

Instead of $CO_2$, in such system can be used other gases (helium, neon, argon, nitrogen, etc.). However, it should be noted that any other gas in comparison to $CO_2$ has also different thermodynamic properties and the observed system will not operate with the same operating parameters (with the same temperatures, pressures and mass flow rates) in each operating point from Fig. 1, if any other mentioned gas is used. One of the possibilities in future research of the presented system will surely be investigation of change in operating parameters and system performances when other gases are used.

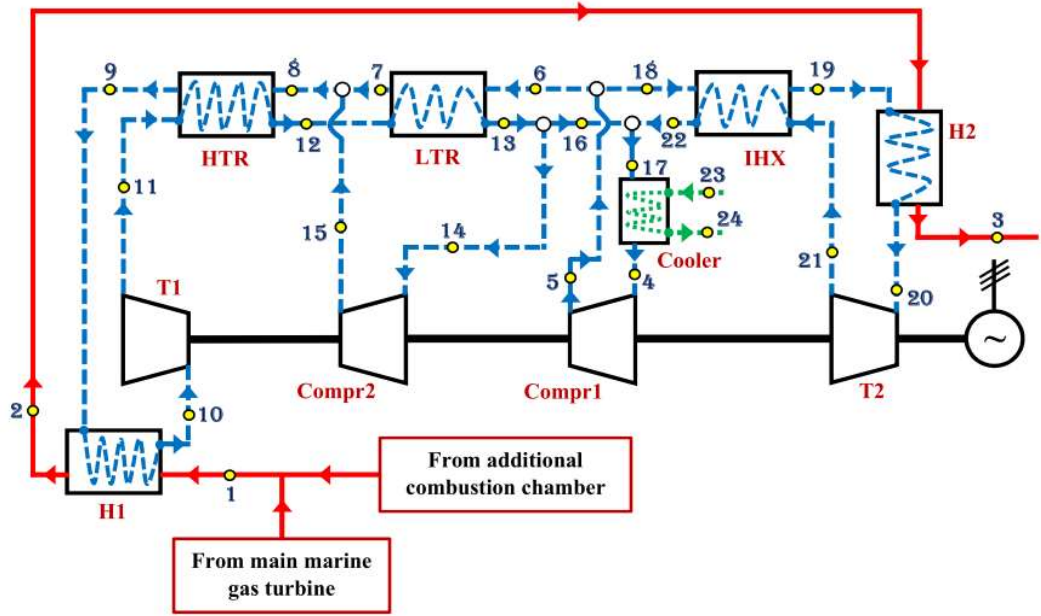

**Fig. 1** Scheme of the analyzed waste heat recovery $CO_2$ closed-cycle gas turbine system along with marked operating points necessary for the exergy analysis
**Source:** Authors

Exergy analysis of the whole observed waste heat recovery $CO_2$ closed-cycle gas turbine system is performed by using scheme presented in Fig. 2. Exergy power inputs into the whole observed system are combustion gases from the gas turbine (and from the additional combustion chamber) and cooling water at the Cooler inlet. Exergy power outputs from the whole observed system are combustion gases at the outlet of the observed system (after performing both heat transfers in H1 and H2, Fig. 1) as well as cooling water at the Cooler outlet and useful produced mechanical power. Those elements are sufficient for the exergy analysis of the whole observed system.

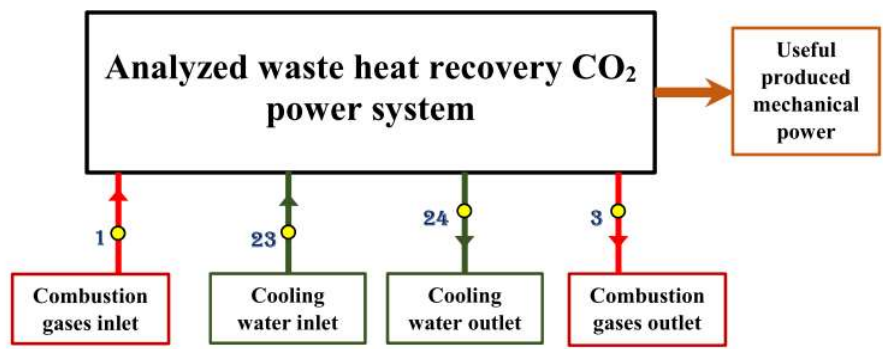

**Fig. 2** Scheme for defining exergy analysis parameters of the whole observed waste heat recovery $CO_2$ closed-cycle gas turbine system
**Source:** Authors

### 3 Operating parameters required for the $CO_2$ closed-cycle gas turbine system exergy analysis

For the exergy analysis of the whole observed waste heat recovery $CO_2$ closed-cycle gas turbine system as well as for the exergy analysis of each its component are required pressures, temperatures and mass flow rates in each operating point from Fig. 1 and Fig. 2. These operating parameters are found in [27] and presented in Table 1.

Specific enthalpies and specific entropies of each fluid stream from Table 1 are calculated by using NIST-REFPROP 9.0 software [28] from known pressures and temperatures.

Specific exergies of each fluid stream are calculated by using Eq. 4. For the specific exergies calculation must be defined the ambient base (dead) state for which the calculations are performed [29, 30]. Ambient base (dead) state for the exergy analysis can be selected provisionally [31, 32] and is related to the pressure and temperature of the ambient in which system or a control volume operates. In this paper, with the ambient pressure of 1 bar (100 kPa) and the ambient temperature of 15 °C (288.15 K) is defined base (dead) state for the exergy analysis.

With an aim to simplify calculations, in this paper in operating points 1, 2 and 3 are used properties of pure air instead of combustion gases properties. This simplification will surely bring some differences in comparison to the case when the combustion gases were used, but the differences will be small because the combustion gases from the main gas turbine as well as combustion gases from the additional combustion chamber, Fig. 1, are dominantly composed of air with a small proportion of fuel [33, 34]. Therefore, this simplification will not have a major influence on the obtained exergy analysis results for any component as well as for the whole observed system, and at the same time will make performed calculations easier and faster.

**Table 1** Fluid stream operating parameters of the analyzed waste heat recovery $CO_2$ closed-cycle gas turbine system [27]

| O.P.* | Medium | Temperature (°C) | Pressure (kPa) | Mass flow rate (kg/s) | Specific enthalpy (kJ/kg) | Specific entropy (kJ/kg·K) | Specific exergy (kJ/kg) |
|---|---|---|---|---|---|---|---|
| 1 | Air | 531.8 | 104.3 | 93.8 | 953.85 | 4.9019 | 236.34 |
| 2 | Air | 369.9 | 104.3 | 93.8 | 778.89 | 4.6594 | 131.26 |
| 3 | Air | 183.1 | 104.3 | 93.8 | 584.51 | 4.3027 | 39.672 |
| 4 | $CO_2$ | 31.0 | 7579 | 127.9 | 298.80 | 1.3226 | 201.33 |
| 5 | $CO_2$ | 60.6 | 21190 | 127.9 | 323.08 | 1.3372 | 221.40 |
| 6 | $CO_2$ | 60.6 | 21190 | 55.7 | 323.08 | 1.3372 | 221.40 |
| 7 | $CO_2$ | 171.5 | 21190 | 55.7 | 552.37 | 1.9383 | 277.47 |
| 8 | $CO_2$ | 171.5 | 21190 | 67.1 | 552.37 | 1.9383 | 277.47 |
| 9 | $CO_2$ | 303.9 | 21190 | 67.1 | 730.09 | 2.2896 | 353.96 |
| 10 | $CO_2$ | 501.8 | 21190 | 67.1 | 974.51 | 2.6539 | 493.42 |
| 11 | $CO_2$ | 384.1 | 7579 | 67.1 | 849.65 | 2.6751 | 362.44 |
| 12 | $CO_2$ | 227.3 | 7579 | 67.1 | 671.93 | 2.3664 | 273.68 |
| 13 | $CO_2$ | 70.6 | 7579 | 67.1 | 481.64 | 1.9062 | 215.99 |
| 14 | $CO_2$ | 70.6 | 7579 | 11.4 | 481.64 | 1.9062 | 215.99 |
| 15 | $CO_2$ | 171.5 | 21190 | 11.4 | 552.37 | 1.9383 | 277.47 |
| 16 | $CO_2$ | 70.6 | 7579 | 55.7 | 481.64 | 1.9062 | 215.99 |
| 17 | $CO_2$ | 70.6 | 7579 | 127.9 | 481.64 | 1.9062 | 215.99 |
| 18 | $CO_2$ | 60.6 | 21190 | 72.2 | 323.08 | 1.3372 | 221.40 |
| 19 | $CO_2$ | 153.0 | 21190 | 72.2 | 523.05 | 1.8709 | 267.56 |
| 20 | $CO_2$ | 339.9 | 21190 | 72.2 | 774.92 | 2.3650 | 377.08 |
| 21 | $CO_2$ | 235.9 | 7579 | 72.2 | 681.59 | 2.3855 | 277.83 |
| 22 | $CO_2$ | 70.6 | 7579 | 72.2 | 481.64 | 1.9062 | 215.99 |
| 23 | Water | 25.0 | 200 | 1119.0 | 105.01 | 0.3672 | 0.8095 |
| 24 | Water | 30.0 | 200 | 1119.0 | 125.91 | 0.4367 | 1.6781 |

* O.P. = Operating Point (according to presented operating points in Fig. 1 and Fig. 2)

## 4 Exergy analysis equations

In comparison to energy analysis which did not take into consideration parameters of the ambient in which analyzed system or a control volume operates [35, 36], exergy analysis takes into consideration parameters of the ambient (ambient pressure and ambient temperature) in which system or a control volume operates [37, 38]. Therefore, exergy analysis takes into consideration additional losses related to the ambient, which are not included in the energy analysis [39]. According to this fact, exergy analysis of any system or a component (control volume) present better operation and losses overview, so it can be a baseline for proper economy (exergo-economy) analysis [40-42].

### 4.1 Overall equations and balances valid for any system or a control volume

As defined in [43, 44], exergy analysis of any system or a control volume has a baseline in the second law of thermodynamics (second law analysis). Overall exergy balance equation, valid for any system or a control volume, can be written according to [45, 46] as:

$$\dot{X}_H + P_{IN} + \sum \dot{E}x_{IN} = P_{OUT} + \sum \dot{E}x_{OUT} + \dot{E}x_D \qquad (1)$$

From Eq. 1 two variables should be defined by an additional equations. The first of these variables is $\dot{X}_H$ – exergy transfer by heat at the temperature $T$, which is defined according to [47, 48] as:

$$\dot{X}_H = \sum (1 - \frac{T_0}{T}) \cdot \dot{Q} \qquad (2)$$

The second variable from Eq. 1 which needs additional definition is $\dot{E}x$ - a total exergy flow of any fluid stream, which equation is, according to [49]:

$$\dot{E}x = \dot{m} \cdot \varepsilon \qquad (3)$$

In Eq. 3, $\varepsilon$ is specific exergy of any fluid stream which is calculated according to recommendations from [50], by using an equation:

$$\varepsilon = (h - h_0) - T_0 \cdot (s - s_0) \qquad (4)$$

The exergy efficiency of any system or a control volume can be defined through general form as recommended in [51, 52], by an equation:

$$\eta x = \frac{\text{cumulative exergy output}}{\text{cumulative exergy input}} \qquad (5)$$

In standard operation of any system or a control volume, mass flow rate leakage of any fluid stream usually did not occur [53]. So, in standard operation, the mass flow rate balance [54] can also be applied:

$$\sum \dot{m}_{IN} = \sum \dot{m}_{OUT} \qquad (6)$$

Presented overall exergy equations and balances will be a baseline for defining all exergy equations during the analysis of whole observed waste heat recovery $CO_2$ closed-cycle gas turbine system as well as for defining exergy analysis equations of each system component.

**4.2 Equations for the exergy analysis of the observed $CO_2$ closed-cycle gas turbine system**

Equations required for the exergy analysis of the whole observed system and for each system component are arranged in several tables in this sub-section. For each component, as well as for the whole observed system are defined equations for the calculation of exergy power input, exergy power output, exergy destruction (exergy power loss) and exergy efficiency, as a standard elements in any exergy analysis [55, 56].
Markings in all the equations presented in this sub-section related to each component are defined according to operating points from Fig. 1, while markings in all the equations presented in this sub-section related to the whole observed system are defined according to operating points from Fig. 2.
In Table 2 are presented exergy analysis equations of all heat exchangers from the observed waste heat recovery $CO_2$ closed-cycle gas turbine system. Equations for each heat exchanger are defined according to recommendations from the literature [57, 58].

**Table 2** Exergy analysis equations of all heat exchangers from the observed waste heat recovery $CO_2$ closed-cycle gas turbine system

| Component | Exergy power input | Eq. | Exergy power output | Eq. |
|---|---|---|---|---|
| H1 | $\dot{E}x_{IN,H1} = \dot{E}x_1 - \dot{E}x_2$ | (7) | $\dot{E}x_{OUT,H1} = \dot{E}x_{10} - \dot{E}x_9$ | (13) |
| H2 | $\dot{E}x_{IN,H2} = \dot{E}x_2 - \dot{E}x_3$ | (8) | $\dot{E}x_{OUT,H2} = \dot{E}x_{20} - \dot{E}x_{19}$ | (14) |
| HTR | $\dot{E}x_{IN,HTR} = \dot{E}x_{11} - \dot{E}x_{12}$ | (9) | $\dot{E}x_{OUT,HTR} = \dot{E}x_9 - \dot{E}x_8$ | (15) |
| LTR | $\dot{E}x_{IN,LTR} = \dot{E}x_{12} - \dot{E}x_{13}$ | (10) | $\dot{E}x_{OUT,LTR} = \dot{E}x_7 - \dot{E}x_6$ | (16) |
| IHX | $\dot{E}x_{IN,IHX} = \dot{E}x_{21} - \dot{E}x_{22}$ | (11) | $\dot{E}x_{OUT,IHX} = \dot{E}x_{19} - \dot{E}x_{18}$ | (17) |
| Cooler | $\dot{E}x_{IN,COOL} = \dot{E}x_{17} - \dot{E}x_4$ | (12) | $\dot{E}x_{OUT,COOL} = \dot{E}x_{24} - \dot{E}x_{23}$ | (18) |
| Component | Exergy destruction | Eq. | Exergy efficiency | Eq. |
| H1 | $\dot{E}x_{D,H1} = \dot{E}x_{IN,H1} - \dot{E}x_{OUT,H1}$ | (19) | $\eta x_{H1} = \dfrac{\dot{E}x_{OUT,H1}}{\dot{E}x_{IN,H1}}$ | (25) |
| H2 | $\dot{E}x_{D,H2} = \dot{E}x_{IN,H2} - \dot{E}x_{OUT,H2}$ | (20) | $\eta x_{H2} = \dfrac{\dot{E}x_{OUT,H2}}{\dot{E}x_{IN,H2}}$ | (26) |
| HTR | $\dot{E}x_{D,HTR} = \dot{E}x_{IN,HTR} - \dot{E}x_{OUT,HTR}$ | (21) | $\eta x_{HTR} = \dfrac{\dot{E}x_{OUT,HTR}}{\dot{E}x_{IN,HTR}}$ | (27) |
| LTR | $\dot{E}x_{D,LTR} = \dot{E}x_{IN,LTR} - \dot{E}x_{OUT,LTR}$ | (22) | $\eta x_{LTR} = \dfrac{\dot{E}x_{OUT,LTR}}{\dot{E}x_{IN,LTR}}$ | (28) |
| IHX | $\dot{E}x_{D,IHX} = \dot{E}x_{IN,IHX} - \dot{E}x_{OUT,IHX}$ | (23) | $\eta x_{IHX} = \dfrac{\dot{E}x_{OUT,IHX}}{\dot{E}x_{IN,IHX}}$ | (29) |
| Cooler | $\dot{E}x_{D,COOL} = \dot{E}x_{IN,COOL} - \dot{E}x_{OUT,COOL}$ | (24) | $\eta x_{COOL} = \dfrac{\dot{E}x_{OUT,COOL}}{\dot{E}x_{IN,COOL}}$ | (30) |

**Source:** Authors

Inside the observed waste heat recovery $CO_2$ closed-cycle gas turbine system exists two $CO_2$ gas turbines (T1 and T2) as well as two $CO_2$ compressors (Compr1 and Compr2). Produced mechanical power by each turbine, as well as used mechanical power by each compressor is an essential element in exergy analysis of each component. Therefore, the equations for the calculation of produced and used mechanical power by each turbine or compressor are defined according to recommendations from [59] and presented in Table 3.

Equations for all exergy analysis parameters of each turbine and each compressor from the observed waste heat recovery $CO_2$ closed-cycle gas turbine system are presented in Table 4 and are obtained by using the instructions from the literature [60].

**Table 3** Equations for mechanical power calculation – produced (turbines), used (compressors) and useful power

| Component | Mechanical power | Eq. |
|---|---|---|
| T1 | $P_{T1} = \dot{m}_{10} \cdot (h_{10} - h_{11})$ | (31) |
| T2 | $P_{T2} = \dot{m}_{20} \cdot (h_{20} - h_{21})$ | (32) |
| Compr1 | $P_{Compr1} = \dot{m}_4 \cdot (h_5 - h_4)$ | (33) |
| Compr2 | $P_{Compr2} = \dot{m}_{14} \cdot (h_{15} - h_{14})$ | (34) |
| Useful | $P_{Useful} = P_{T1} + P_{T2} - P_{Compr1} - P_{Compr2}$ | (35) |

**Source:** Authors

**Table 4** Exergy analysis equations of all turbines and compressors from the observed waste heat recovery $CO_2$ closed-cycle gas turbine system

| Component | Exergy power input | Eq. | Exergy power output | Eq. |
|---|---|---|---|---|
| T1 | $\dot{E}x_{IN,T1} = \dot{E}x_{10}$ | (36) | $\dot{E}x_{OUT,T1} = \dot{E}x_{11} + P_{T1}$ | (40) |
| T2 | $\dot{E}x_{IN,T2} = \dot{E}x_{20}$ | (37) | $\dot{E}x_{OUT,T2} = \dot{E}x_{21} + P_{T2}$ | (41) |
| Compr1 | $\dot{E}x_{IN,Compr1} = \dot{E}x_4 + P_{Compr1}$ | (38) | $\dot{E}x_{OUT,Compr1} = \dot{E}x_5$ | (42) |
| Compr2 | $\dot{E}x_{IN,Compr2} = \dot{E}x_{14} + P_{Compr2}$ | (39) | $\dot{E}x_{OUT,Compr2} = \dot{E}x_{15}$ | (43) |
| **Component** | **Exergy destruction** | **Eq.** | **Exergy efficiency** | **Eq.** |
| T1 | $\dot{E}x_{D,T1} = \dot{E}x_{IN,T1} - \dot{E}x_{OUT,T1}$ | (44) | $\eta x_{T1} = \dfrac{P_{T1}}{\dot{E}x_{10} - \dot{E}x_{11}}$ | (48) |
| T2 | $\dot{E}x_{D,T2} = \dot{E}x_{IN,T2} - \dot{E}x_{OUT,T2}$ | (45) | $\eta x_{T2} = \dfrac{P_{T2}}{\dot{E}x_{20} - \dot{E}x_{21}}$ | (49) |
| Compr1 | $\dot{E}x_{D,Compr1} = \dot{E}x_{IN,Compr1} - \dot{E}x_{OUT,Compr1}$ | (46) | $\eta x_{Compr1} = \dfrac{\dot{E}x_5 - \dot{E}x_4}{P_{Compr1}}$ | (50) |
| Compr2 | $\dot{E}x_{D,Compr2} = \dot{E}x_{IN,Compr2} - \dot{E}x_{OUT,Compr2}$ | (47) | $\eta x_{Compr2} = \dfrac{\dot{E}x_{15} - \dot{E}x_{14}}{P_{Compr2}}$ | (51) |

**Source:** Authors

In Table 5 are presented all the required exergy analysis equations related to the whole observed waste heat recovery $CO_2$ closed-cycle gas turbine system (markings are related to operating points from Fig. 2). As for each component, in any equation from Table 5 related to the whole system, are followed recommendations from the literature [61, 62].

**Table 5** Equations for the exergy analysis of the whole observed waste heat recovery $CO_2$ closed-cycle gas turbine system

| - - - | Whole system | Eq. |
|---|---|---|
| **Exergy power input** | $\dot{E}x_{IN,WS} = \dot{E}x_1 + \dot{E}x_{23}$ | (52) |
| **Exergy power output** | $\dot{E}x_{OUT,WS} = \dot{E}x_3 + \dot{E}x_{24} + P_{Useful}$ | (53) |
| | $\dot{E}x_{D,WS} = \dot{E}x_{IN,WS} - \dot{E}x_{OUT,WS}$ | (54) |
| **Exergy destruction** | $\dot{E}x_{D,WS} = \sum \dot{E}x_{D,\text{of all components}} =$ $\dot{E}x_{D,H1} + \dot{E}x_{D,H2} + \dot{E}x_{D,HTR} + \dot{E}x_{D,LTR} +$ $\dot{E}x_{D,IHX} + \dot{E}x_{D,COOL} + \dot{E}x_{D,T1} + \dot{E}x_{D,T2} +$ $\dot{E}x_{D,Compr1} + \dot{E}x_{D,Compr2}$ | (55) |
| **Exergy efficiency** | $\eta x_{WS} = \dfrac{P_{Useful}}{\dot{E}x_1 + \dot{E}x_{23} - \dot{E}x_3 - \dot{E}x_{24}}$ | (56) |

**Source:** Authors

## 5 Results and discussion

Exergy power inputs and outputs of all heat exchangers from the observed waste heat recovery $CO_2$ closed-cycle gas turbine system are presented in Fig. 3. Exact values of each exergy power input and output for each heat exchanger are presented under the diagram in Fig. 3.

It should be highlighted that H1 has the highest exergy power input and output in comparison to other heat exchangers because it is the first heat exchanger to whom are delivered combustion gases from the main marine gas turbine and from the additional combustion chamber, Fig. 1. Heat exchanger H2 has a little lower exergy power input and output in comparison to H1, Fig. 3. However, from Fig. 3 can be concluded that heat exchangers in which heat is transferred from combustion gases to $CO_2$ have the highest exergy power inputs and outputs in the observed system.

Cooler, used for $CO_2$ cooling before its compression in Compr1 ($CO_2$ cooling is performed with water) has the lowest exergy power input and output in comparison to all other heat exchangers in the observed system.

In this analysis, presentation of exergy power inputs and outputs for all heat exchangers from the analyzed waste heat recovery $CO_2$ closed-cycle gas turbine system is crucial, because the difference between exergy power input and output define exergy destruction (exergy power loss), while the ratio of exergy power output and input define exergy efficiency of each heat exchanger, Table 2.

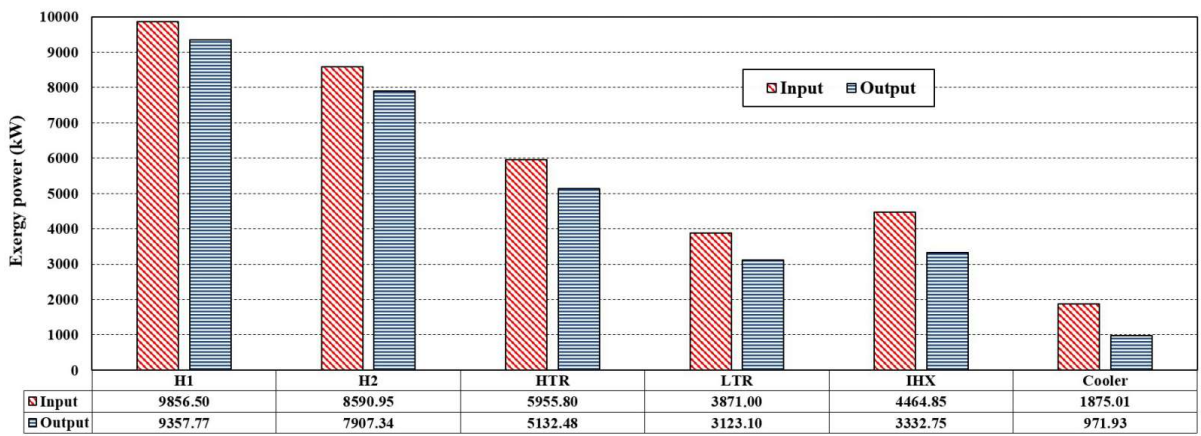

**Fig. 3** Exergy power inputs and outputs of all heat exchangers from the observed waste heat recovery $CO_2$ closed-cycle gas turbine system

**Source:** Authors

Exergy destructions (exergy power losses) and exergy efficiencies of all heat exchangers in the observed waste heat recovery $CO_2$ closed-cycle gas turbine system are presented in Fig. 4.

Heat exchangers H1 and H2 in which the heat is transferred from combustion gases to $CO_2$ have the lowest exergy destructions equal to 498.74 kW for H1 and 683.61 kW for H2 in comparison to all other heat exchangers. Consequentially, H1 and H2 have the highest exergy efficiencies in comparison to all other heat exchangers (94.94% for H1 and 92.04% for H2). Heat exchangers in which hotter $CO_2$ transferred heat to colder $CO_2$ (HTR, LTR and IHX) have higher exergy destructions and consequentially lower exergy efficiencies in comparison to H1 and H2, Fig. 4. It is interesting to observe that HTR has higher exergy efficiency equal to 86.18% in comparison to LTR (80.68%), regardless of higher exergy destruction (823.32 kW for HTR and 747.90 kW for LTR). The highest exergy destruction in comparison to other heat exchangers from the observed system is detected in IHX and is equal to 1132.10 kW, Fig. 4. Also, it should be noted that regardless of the highest exergy destruction, IHX did not have the lowest exergy efficiency when taking into account all heat exchangers from the observed system.

When taking into consideration all the heat exchangers from the analyzed waste heat recovery $CO_2$ closed-cycle gas turbine system, the lowest exergy efficiency is detected in Cooler and equal to 51.84%. The Cooler has exergy destruction lower than IHX, but higher in comparison to other heat exchangers from the observed system, so the Cooler exergy destruction cannot be the only reason for such low exergy efficiency. The main reason for low Cooler exergy efficiency is cooling water, which represents Cooler exergy power output, Table 2. Water in the Cooler has significantly lower exergy power in comparison to $CO_2$ (Table 1) what resulted with a low value of Cooler exergy efficiency.

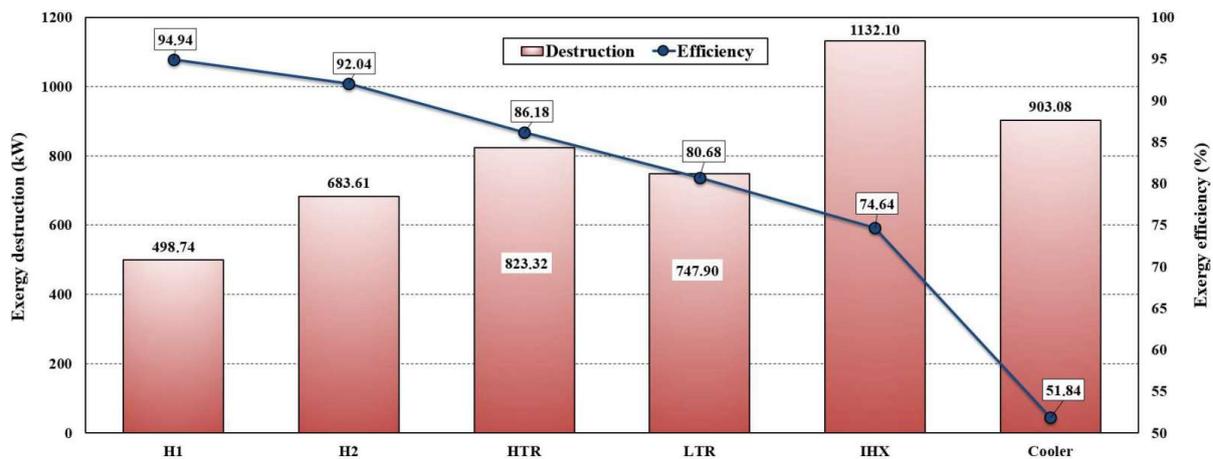

**Fig. 4** Exergy destruction and exergy efficiency of all heat exchangers from the observed waste heat recovery $CO_2$ closed-cycle gas turbine system
**Source:** Authors

In the analyzed waste heat recovery $CO_2$ closed-cycle gas turbine system, two $CO_2$ gas turbines (T1 and T2) are mechanical power producers, while two $CO_2$ compressors (Compr1 and Compr2) are mechanical power consumers. The difference in produced mechanical power by turbines and consumed mechanical power by compressors is delivered to an electrical generator for additional electrical power production. Without the observed waste heat recovery system, heat energy from combustion gases will be released to the atmosphere and will be completely lost.

Fig. 5 presents produced, used and useful mechanical power in the waste heat recovery $CO_2$ closed-cycle gas turbine system. Each mechanical power is calculated by using equations from Table 3.

From the viewpoint of gas turbines as the mechanical power producers, it can be concluded that the dominant mechanical power producer is T1, while the secondary gas turbine is T2 (T1 produces 8378.11 kW, while T2 produces 6738.43 kW of mechanical power), Fig. 5. This difference can be explained by using Table 1 – T1 operates with little lower $CO_2$ mass flow rate, but with significantly higher $CO_2$ temperature at the inlet in comparison to T2.

The dominant mechanical power consumer is Compr1 which consumes 3105.41 kW of mechanical power, while Compr2 consumes only 806.32 kW of mechanical power. This difference can be explained again by using Table 1 – Compr1 operates with significantly higher $CO_2$ mass flow rate in comparison to Compr2.

The difference in produced mechanical power (cumulative by both gas turbines) and in used mechanical power (cumulative by both compressors) is useful power, which for the observed system and according to $CO_2$ operating parameters from Table 1, equals 11204.80 kW. Useful mechanical power drives an electrical generator in the analyzed system.

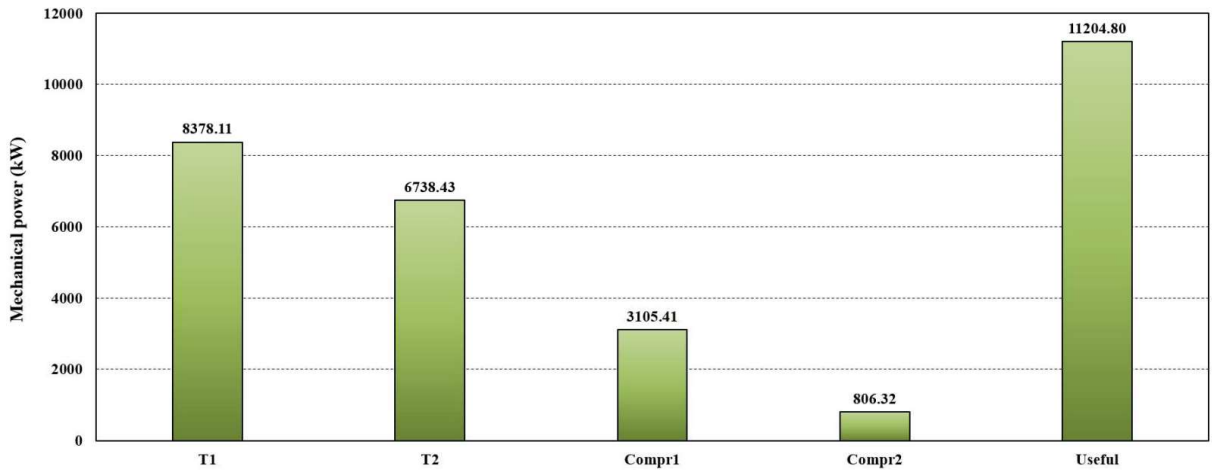

**Fig. 5** Produced (turbines), used (compressors) and useful mechanical power from the observed waste heat recovery $CO_2$ closed-cycle gas turbine system
**Source:** Authors

Exergy power inputs and outputs for all turbines and compressors from the observed system are calculated by using equations from Table 4. Necessary part of turbines and compressors exergy power inputs and outputs are also produced or used mechanical power presented in Fig. 5. For each turbine and compressor from the observed system, exergy power inputs and outputs are presented in Fig. 6, while the exact values of each input and output are presented under the diagram of each component.

From Fig. 6 can be clearly seen that the highest exergy power inputs and outputs of all turbines and compressors have T1, while the lowest exergy power inputs and outputs is observed for Compr2.

A comparison of both turbines (T1 and T2) shows that T1 has higher exergy power inputs and outputs, but the difference is not significant. A comparison of both compressors (Compr1 and Compr2) shows that Compr1 has significantly higher exergy power inputs and outputs, what can be explained with several times higher $CO_2$ mass flow rate through Compr1, Table 1. Also, it is interesting to observe that Compr1 has higher exergy power input and output in comparison to T2, Fig. 6 – the main reason of such occurrence is almost two times higher $CO_2$ mass flow rate through Compr1 than through T2, regardless of the fact that T2 operates with much higher $CO_2$ temperatures, Table 1.

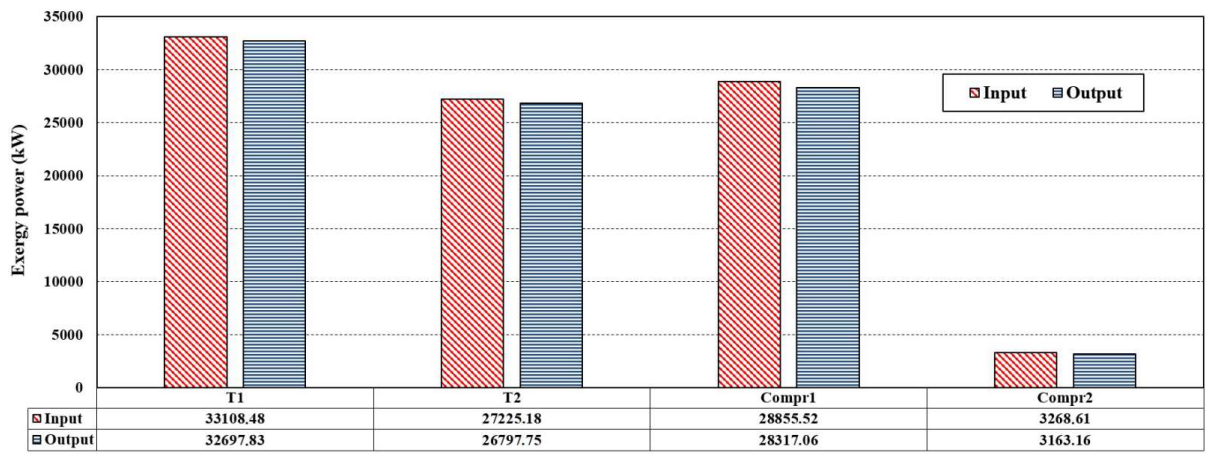

**Fig. 6** Exergy power inputs and outputs of all turbines and compressors from the observed waste heat recovery $CO_2$ closed-cycle gas turbine system
**Source:** Authors

Exergy destructions and exergy efficiencies of all turbines and compressors from the analyzed waste heat recovery $CO_2$ closed-cycle gas turbine system are calculated by using equations from Table 4 and presented in Fig. 7.

By observing turbines from the analyzed system, T2 has higher exergy destruction and consequentially lower exergy efficiency in comparison to T1 (exergy destructions of T1 and T2 are equal to 410.65 kW and 427.42 kW, while exergy efficiencies are equal to 95.33% and 94.04%, respectively). However, it should be highlighted that the differences in exergy destructions and in exergy efficiencies of both turbines are small.

Comparison of Compr1 and Compr2 shows that Compr1 has notably higher exergy destruction and consequentially notably lower exergy efficiency. Therefore the differences in exergy destruction and in exergy efficiency between observed compressors are much higher in comparison to the same differences between turbines. Compr1 has exergy destruction and exergy efficiency equal to 538.46 kW and 82.66%, while the same parameters for Compr2 are 105.45 kW and 86.92%, respectively, Fig. 7.

Comparison of exergy efficiencies between mechanical power producers (turbines) and mechanical power consumers (compressors), Fig. 7, show that in the observed system both turbines have notably higher exergy efficiencies in comparison to any compressor.

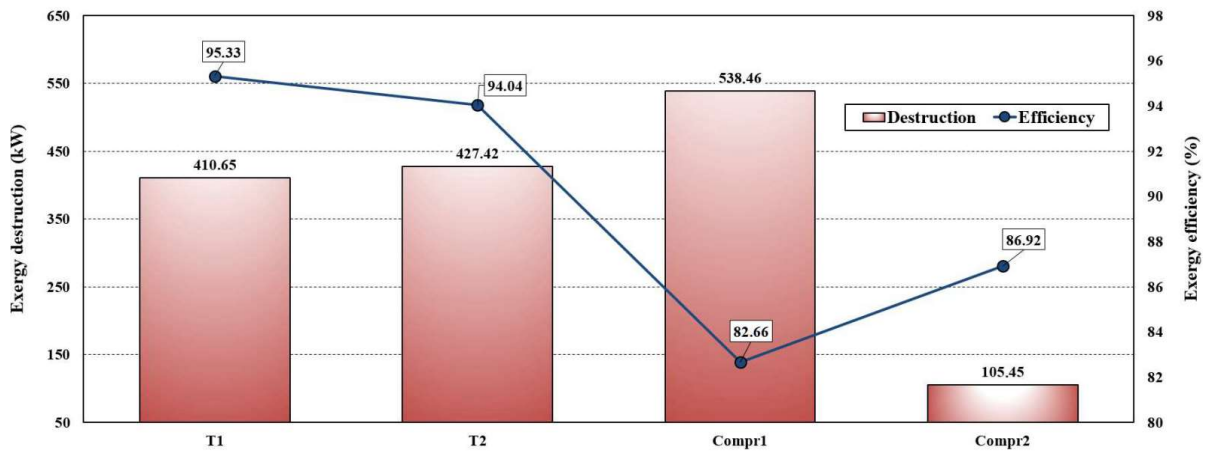

**Fig. 7** Exergy destructions and exergy efficiencies of all turbines and compressors from the observed waste heat recovery $CO_2$ closed-cycle gas turbine system
**Source:** Authors

Exergy analysis of the whole observed waste heat recovery $CO_2$ closed-cycle gas turbine system is performed according to the scheme presented in Fig. 2 and by using equations from Table 5. Obtained exergy analysis results of the whole observed system, by using operating parameters from Table 1, are presented in Table 6.

From Table 6 can clearly be seen that the whole observed waste heat recovery $CO_2$ closed-cycle gas turbine system has exergy destruction equal to 6270.73 kW, while the exergy efficiency is equal to 64.12%.

Obtained exergy efficiency of the whole observed system is higher in comparison to conventional closed-cycle gas turbines [63, 64], what proves that application of such waste heat recovery system on-board ships with the gas turbine as a main propulsion element will be beneficial for the additional electrical power production.

However, it should be highlighted that such system will also bring many disadvantages into the ship engine room. The whole system uses $CO_2$ as an operating medium, what will require installing of additional $CO_2$ tanks and compressors for $CO_2$ delivery to the system and $CO_2$ removing from the system (useful produced mechanical power in such systems is regulated by change of operating medium mass flow rate). Such system will require additional space in the ship engine room. As the operating medium is high pressure and high temperature $CO_2$, it should be ensured a good sealing in all the places where $CO_2$ can exit from the system. Maintenance of such system is complex and will require additional time and effort of the ship crew. Regardless of all the disadvantages, the additional produced mechanical power can overcome all the disadvantages, especially in a long operation period.

Presented waste heat recovery system will surely increase the efficiency of the whole propulsion plant with the gas turbine as main propulsion element and will ensure competitiveness of such propulsion plant with other known marine propulsion systems.

**Table 6** Obtained exergy analysis parameters of the whole observed waste heat recovery $CO_2$ closed-cycle gas turbine system

| Whole system | Exergy power input (kW) | Exergy power output (kW) | Exergy destruction (kW) | Exergy efficiency (%) |
|---|---|---|---|---|
| | 23074.56 | 16803.83 | 6270.73 | 64.12 |

**Source:** Authors

Future research of the analyzed waste heat recovery $CO_2$ closed-cycle gas turbine system will be based on the application of AI (Artificial Intelligence) methods and processes already developed by our research team [65-68]. The intention will be to analyze and perform optimization of the system operating parameters (Table 1) with an aim to decrease the system exergy destruction and simultaneously to increase system exergy efficiency.

## 6 Conclusions

In this paper is performed an exergy analysis of marine waste heat recovery $CO_2$ closed-cycle gas turbine system. The presented system operates by using combustion gases from the main marine gas turbine and from the additional combustion chamber (if required). It is performed analysis of each system component individually, as well as analysis of the whole observed system. It is detected components which operation can be improved. The main conclusions from the presented research are:
- <u>Heat exchangers:</u>  Heat exchangers H1 and H2 in comparison to all other heat exchangers have the lowest exergy destructions equal to 498.74 kW for H1 and 683.61 kW for H2 and simultaneously the highest exergy efficiencies (94.94% for H1 and 92.04% for H2). Heat exchangers HTR, LTR and IHX have higher exergy destructions and consequentially lower exergy efficiencies in comparison to H1 and H2. The lowest exergy efficiency is detected in Cooler (51.84%) because cooling water has significantly lower exergy power in comparison to $CO_2$.
- <u>Turbines and Compressors:</u> The dominant mechanical power producer is T1, while the secondary gas turbine is T2 (T1 produces 8378.11 kW, while T2 produces 6738.43 kW of mechanical power). T2 has higher exergy destruction and lower exergy efficiency in comparison to T1 (exergy destructions of T1 and T2 are equal to 410.65 kW and 427.42 kW, while exergy efficiencies are equal to 95.33% and 94.04%, respectively). The dominant mechanical power consumer is Compr1 which consumes 3105.41 kW of mechanical power, while Compr2 consumes only 806.32 kW of mechanical power. Compr1 has notably higher exergy destruction and notably lower exergy efficiency in comparison to Compr2. Exergy destruction and exergy efficiency of Compr1 are equal to 538.46 kW and 82.66%, while the same parameters for Compr2 are 105.45 kW and 86.92%, respectively.
- <u>Whole waste heat recovery system:</u> The whole observed waste heat recovery system has exergy destruction equal to 6270.73 kW, while the exergy efficiency is 64.12%. Useful power produced by the whole system and used for electrical generator drive equals 11204.80 kW.

Obtained high exergy efficiency of the whole observed system proves its application on-board ships with the gas turbine as a main propulsion element, for the additional electrical power production.

**Acknowledgment**

This research has been supported by the Croatian Science Foundation under the project IP-2018-01-3739, CEEPUS network CIII-HR-0108, European Regional Development Fund under the grant KK.01.1.1.01.0009 (DATACROSS), project CEKOM under the grant KK.01.2.2.03.0004, CEI project "COVIDAi" (305.6019-20), University of Rijeka scientific grant uniri-tehnic-18-275-1447 and University of Rijeka scientific grant uniri-tehnic-18-18-1146.

| **NOMENCLATURE** | |
|---|---|
| | |
| *Latin Symbols:* | |
| $\dot{E}x$ | total exergy flow of fluid, kW |
| $h$ | specific enthalpy, kJ/kg |
| $\dot{m}$ | mass flow rate, kg/s |
| $p$ | pressure, bar or kPa |
| $P$ | mechanical power, kW |
| $\dot{Q}$ | energy transfer by heat, kW |
| $s$ | specific entropy, kJ/kg·K |
| $T$ | temperature, °C or K |
| $\dot{X}_H$ | exergy transfer by heat, kW |
| | |
| *Greek symbols:* | |
| $\varepsilon$ | specific exergy, kJ/kg |
| $\eta x$ | exergy efficiency, - |
| | |

| Subscripts: | |
|---|---|
| 0 | base ambient state (dead state) |
| D | exergy destruction (exergy loss) |
| IN | inlet (input) |
| OUT | outlet (output) |
| WS | whole system |
| | |
| Abbreviations: | |
| Compr1 | $CO_2$ compressor 1 |
| Compr2 | $CO_2$ compressor 2 |
| COOL | water cooler |
| H1 | combustion gases heater 1 |
| H2 | combustion gases heater 2 |
| HTR | high temperature recuperator |
| IHX | intermediate heat exchanger |
| LTR | low temperature recuperator |
| T1 | $CO_2$ turbine 1 |
| T2 | $CO_2$ turbine 2 |